\documentclass{amsart}
\usepackage{amssymb,amsmath,amsthm,xspace}%\usepackage{showkeys}
\usepackage[small]{diagrams}
\def\mat#1{\ensuremath{#1}\xspace}

\def\cF{\mat{\mathbb{F}}}
\def\cN{\mat{\mathbb{N}}}   %natural numbers
\def\cQ{\mat{\mathbb{Q}}}   %ring of rationals
\def\cC{\mat{\mathbb{C}}}   %complex numbers

\def\cZ{\mat{\mathbb{Z}}}   %ring of integers

\def\lE{\mat{\mathcal{E}}}

\def\lH{\mat{\mathcal{H}}}

\def\lP{\mat{\mathcal{P}}}

\let\Phitemp\Phi \def\Phi{\mat{\Phitemp}}
\let\Psitemp\Psi \def\Psi{\mat{\Psitemp}}
\let\etatemp\eta \def\eta{\mat{\etatemp}}
\def\La{\mat{\Lambda}}
\def\la{\mat{\lambda}}
\def\vi{\mat{\varphi}}
\let\mutemp\mu \def\mu{\mat{\mutemp}}
\let\nutemp\nu \def\nu{\mat{\nutemp}}
\def\si{\mat{\sigma}}

\def\al{\mat{\alpha}}
\def\be{\mat{\beta}}
\def\ga{\mat{\gamma}}
\def\Ga{\mat{\Gamma}}

\def\De{\mat{\Delta}}
\def\Nabla{\bigtriangledown}
\def\Na{\mat{\Nabla}}

\let\xitemp\xi \def\xi{\mat{\xitemp}}

\def\mrm@#1{\mat{\mathrm{#1}}}

\def\g {\mat{\mathfrak{g}}}  %\gg and \ggg are standartly reserved for ">>" and ">>>"
\def\gh{\mat{\mathfrak{h}}}

\def\DMO{\DeclareMathOperator}

\DMO{\Hom}{Hom}
\DMO{\lHom}{\lH\mathit{om}}
\DMO{\Ext}{Ext}
\DMO{\lExt}{\lE\mathit{xt}}
\DMO{\End}{End}
\DMO{\Aut}{Aut}
\DMO{\Fun}{Fun}
\DMO{\Tor}{Tor}
\DMO{\ext}{ext}

\DMO{\Ob}{Ob}
\DMO{\Mor}{Mor}
\DMO{\im}{im}
\DMO{\coim}{coim}
\DMO{\coker}{coker}
\DMO{\Arr}{Arr}

\DMO{\Id}{Id}
\DMO{\add}{add} % splitting of idempotents (karoubinization)
\DMO{\ind}{ind} % category of ind-objects
\DMO{\pro}{pro} % category of pro-objects
\DMO{\Map}{Map} %
\DMO{\Iso}{Iso} %
\DMO{\Isom}{Isom}%
\DMO{\Ind}{Ind}

\DMO{\Presh}{Presh}

\DMO\coalg{Coalg}
\DMO{\Rep}{Rep}
\DMO{\Cor}{Cor}

\DMO{\Mod}{Mod}
\DMO{\rad}{rad}
\DMO{\soc}{soc}
\DMO{\ann}{ann}

\DMO{\Spec}{Spec}
\DMO{\spec}{Spec}
\DMO{\Proj}{Proj}
\DMO{\supp}{supp}
\DMO{\Coh}{Coh}
\DMO{\coh}{Coh}
\DMO{\Qcoh}{QCoh}
\DMO{\QCoh}{QCoh}
\DMO{\Pic}{Pic}
\DMO{\Div}{Div}
\DMO{\ch}{ch}
\DMO{\Hilb}{Hilb}
\DMO{\Fitt}{Fitt}
\DMO{\Quot}{Quot}

\DMO{\Gras}{Gr}
\DMO{\Flag}{Flag}

\DMO{\cone}{cone}
\DMO{\Tw}{Tw}
\DMO{\rank}{rk}
\DMO{\rk}{rk}
\DMO{\codim}{codim}
\DMO{\cov}{cov}
\DMO{\sgn}{sgn}
\DMO{\td}{td}
\DMO{\GL}{GL}
\DMO{\SL}{SL}

\DMO\Der{Der}
\DMO\der{Der}
\DMO\coder{Coder}
\DMO{\diag}{diag}
\DMO{\HMod}{HMod} %the homotopy category of modules over DGC
\DMO{\ad}{ad}
\DMO*{\colim}{colim}
\DMO*{\hocolim}{hocolim}
\DMO*{\holim}{holim}
\DMO{\Ho}{Ho}
\DMO{\har}{char}
\DMO{\sk}{sk}
\DMO{\cosk}{cosk}
\DMO{\Gal}{Gal}
\DMO{\tr}{tr}
\DMO{\Tr}{Tr}
\DMO{\Sh}{Sh}
\DMO{\Is}{Is} %Isometries
\DMO{\Hol}{Hol} %Holomorphic automorphisms
\DMO{\Lie}{Lie} %Lie algebra of a group
\DMO{\Res}{Res} %restriction
\DMO{\irr}{irr} %
\DMO{\Irr}{Irr} %
\DMO{\Exp}{Exp} %
\DMO{\Log}{Log} %
\DMO{\Pow}{Pow}
\DMO{\mult}{mult} %
\DMO{\height}{ht} %
\DMO{\wt}{wt}
\DMO{\Vect}{Vect}
\DMO{\moda}{mod}

\def\dd{\mat{\partial}}

\def\sb{\subset}

\def\xx{\times}

\def\pser#1{[\![#1]\!]} %formal power series [[#1]]
\def\inv{^{-1}}

\def\ub#1{\mat{\overline{#1}}}  %upper bar
\def\set#1{\mat{\{ #1\}}}
\def\sets#1#2{\mat{\{ #1 \mid #2\}}}

\def\arrowsUsual{
\newarrow{TeXto}----{->}
\newarrow{TeXinto}C---{->}
\newarrow{TeXonto}----{->>}
\def\ar{\rightarrow}
\def\emb{\hookrightarrow}
\def\mto{\mapsto}
\def\arr{rTeXto}
\def\embb{\rTeXinto}
\newarrow{Eq}=====
    }

\newif\ifukr\ukrfalse
\newif\ifrus\rusfalse
\newif\ifger\gerfalse

\def\theoremsElsevier{
\newtheorem{nthr}{Auxiliary} %auxiliary theorem, main counter  

\newtheorem{thm}    [nthr]{Theorem}
\newtheorem{lem}    [nthr]{Lemma}
\newtheorem{cor}    [nthr]{Corollary}
\newtheorem{conj}    [nthr]{Conjecture}
\theoremstyle{definition}
\newtheorem{defn}   [nthr]{Definition}
\newtheorem{rem}    [nthr]{Remark}
\newtheorem{exmp}   [nthr]{Example}
}
 
\theoremsElsevier\arrowsUsual

\begin{document}
\title[]{A computational criterion for the Kac conjecture}%

\author{Sergey Mozgovoy}%

%\address{Institut f\"ur Mathematik, Johannes Gutenberg-Universit\"at Mainz,
%55099 Mainz, Germany.}%

\email{mozgov@math.uni-wuppertal.de}%

\thanks{}%
%\subjclass[2000]{16G20}%

%\date{9.08.2006}%
%\dedicatory{}%
%\commby{}%
% ----------------------------------------------------------------
\begin{abstract}
We give a criterion for the Kac conjecture asserting that the free
term of the polynomial counting the absolutely indecomposable
representations of a quiver over a finite field of given dimension
coincides with the corresponding root multiplicity of the
associated Kac-Moody algebra. Our criterion suits very well for
computer tests.
\end{abstract}
\maketitle

\section{Introduction}
Let $\Ga$ be a quiver without loops and let $n$ be its number of
vertices. For any $\al\in\cN^n$, let $m_\al(q)$ and $a_\al(q)$ be
respectively the number of representations and absolutely
indecomposable representations of $\Ga$ over $\cF_q$ of dimension
\al. The underlying graph of $\Ga$ defines a symmetric generalized
$n\xx n$ Cartan matrix $C$ with
$$c_{ij}=
\begin{cases}
2&\text{if }i=j,\\
-b_{ij}&\text{if }i\ne j,
\end{cases}$$
where $b_{ij}$ is the number of edges connecting the vertices $i$
and $j$. Let \g be the corresponding Kac-Moody algebra
\cite{Kac3}. It was shown by Kac \cite{Kac1} that $m_\al(q)$ and
$a_\al(q)$ are polynomials in $q$ with integer coefficients.
Moreover, $a_\al$ is nonzero if and only if \al is a root of \g and
$a_\al=1$ if and only if $\al$ is a real root. It was
conjectured by Kac~\cite{Kac1} that

%{\renewcommand{\thenthr}{1}
\begin{conj}
$a_\al(0)$ is equal to the multiplicity of \al in \g.
\end{conj}
\begin{conj}
All coefficients of $a_\al$ are nonnegative.
\end{conj}

These conjectures were proved by Crawley-Boevey and Van den Bergh
\cite{CBB} in the case when \al is indivisible. In the general
case, there is a criterion for the first Kac conjecture given by
Sevenhant, Van den Bergh and Hua \cite{SB1,Hua1} which, however,
uses the intrinsic structure of the Kac-Moody algebra and is not
very suitable for testing the Kac conjecture on computer. In
contrast to the latter we suggest a criterion which is well
adapted for computer tests.

For any $\al\in\cZ^n$, define $\height\al$ to be the sum of the
coordinates of \al. Define the \cZ-valued quadratic form $T$,
called the Tits form, on $\cZ^n$ by the matrix $\frac 12C$. Let
\lP be the set of partitions. For any multipartition
$\la=(\la^1,\dots,\la^n)\in\lP^n$, we define
$|\la|:=(|\la^1|,\dots,|\la^n|)\in\cN^n$ and
$\la_k:=(\la_k^1,\dots,\la_k^n)\in\cN^n$ for $k\ge1$. Define
$r_\la\in\cQ(q)$ by
$$r_\la(q):=\prod_{k\ge1}\frac{q^{-T(\la_k)}}{\prod_{i=1}^n\vi_{\la^i_k-\la^i_{k+1}}(q\inv)},$$
where $\vi_m(q):=\prod_{i=1}^m(1-q^i)$ for $m\in\cN$. Finally, for
any $\al\in\cN^n$, define $r_\al:=\sum_{|\la|=\al}r_\la$. A
different description of the functions $r_\al\in\cQ(q)$ can be
found in \cite[Section 2]{SB1}. Our main result is the following

\begin{thm}\label{thr intro criterion}
The first Kac conjecture is true if and only if, for any
$\al\in\cN^n$, either $r_\al(0)=0$ or $T(\al)=\height\al$.
\end{thm}

Our proof is based on the formulas relating the functions
$r_\al$ with the polynomials $m_\al$ and $a_\al$ together
with the Peterson recursive formula for the root multiplicities of
the Kac-Moody algebra \cite{Pet1}.

Denote by $\Phi$ the set of all irreducible polynomials in
$\cF_q[t]$ with the leading coefficient $1$, excluding $t$. Denote
by $\Phi_d(q)$ the number of polynomials in $\Phi$ having degree $d$;
it is a polynomial in $q$. The next formula, relating $r_\la$
and $m_\al$, is due to Kac and Stanley \cite[p. 90]{Kac1}
$$m_\al(q)=\sum_{\substack{\nu:\Phi\to\lP^n\\ |\nu|=\al}}\prod_{f\in\Phi}
\psi_{\deg f}(r_{\nu(f)})(q),$$ where $|\nu|:=\sum_{f\in\Phi}\deg
f\cdot|\nu(f)|$ and $\psi_d$ are the Adams operations on $\cQ(q)$ (see
Appendix) given by $\psi_d(f(q))=f(q^d)$ for $f\in\cQ(q)$. 
This formula can be simplified using generating
functions. Consider the functions
$$m(q):=\sum_{\al\in\cN^n}m_\al(q)x^\al,\quad
a(q):=\sum_{\al\in\cN^n}a_\al(q)x^\al,\quad
r(q):=\sum_{\al\in\cN^n}r_\al(q)x^\al$$ %
in $\cQ(q)\pser{x_1,\dots,x_n}$. Then the Kac-Stanley formula can be
written in the form (see Lemma \ref{lmm Kac-St reformulation})
$$m=\prod_{d\ge 1}\psi_d(r)^{\Phi_d}.$$
The functions $m$ and $a$ are related by the formula (see 
Lemma \ref{lmm:rel a and m} and Appendix)  
$$m=\Exp(a).$$
We give a rather technical proof of this formula in Section \ref{sec:relations},
although it has an intuitive explanation. Namely, assume for
simplicity that for all $\al\in\cN^n$ there exists the moduli
space $M_\al$ (respectively, $A_\al$) of representations
(respectively, absolutely indecomposable representations) of
dimension \al and that they have cellular decompositions. (Note
that in general these moduli spaces do not exist). Then the
numbers of cells of fixed dimension correspond to the coefficients
of the polynomials $m_\al$ (respectively, $a_\al$). The fact that
any representation can be uniquely (up to the permutation of
summands) written as a direct sum of indecomposable
representations implies
$$\sum_{\al\in\cN^n} \#M_\al(\cF_q) x^\al
=\prod_{\al\in\cN^n}(1+\#S^1A_\al(\cF_q)x^\al+\#S^2A_\al(\cF_q)x^{2\al}+\dots),$$%
where $S^n$ denotes the symmetric product. Using the existence of
a cellular decomposition of $A_\al$ we obtain
$\#S^nA_\al(\cF_q)=\si_n(a_\al)(q)$. It follows
\begin{align*}
\sum_{\al\in\cN^n} m_\al x^\al
&=\prod_{\al\in\cN^n}(1+\si_1(a_\al)x^\al+\si_2(a_\al)x^{2\al}+\dots)\\
&=\prod_{\al\in\cN^n}\Exp(a_\al x^\al)=\Exp(a).
\end{align*}

The above formulas can be already used for the explicit
calculation of the polynomials $a_\al$. There is, however, a direct relation
between $a$ and $r$ due to Hua (cf. \cite[Theorem 4.6]{Hua1})
$$a(q)=(q-1)\Log(r(q)).$$
Equivalently, this formula can be written as $m(q)=\Pow(r(q),q-1)$ and we show in
Theorem \ref{thr intro relations} that it follows easily from the
Kac-Stanley formula.

In Section \ref{sec:relations} we discuss various relations
between the generating functions $m$, $a$ and $r$. In Section \ref{sec:criterium}
we prove the criterion for the first Kac conjecture using these
relations together with a Peterson recursive formula for the
multiplicities of the Kac-Moody algebra \cite{Pet1}. As the paper
of Peterson \cite{Pet1} is unpublished, we include a rather
detailed description of his approach. In the Appendix we gather
basic definitions concerning the \la-rings. All the computations in the paper
were performed using the algebraic combinatorics package ``MuPAD-Combinat''.

After this paper was finished I became aware of the preprint
\cite{Haus1}, where the proof of the first Kac conjecture is
announced.
\section{Relations between $m$, $a$ and $r$}\label{sec:relations}
The aim of this section is to prove various relations between the
generating functions $m$, $a$ and $r$. As a point of departure,
we use the Kac-Stanley formula
$$m_\al(q)=\sum_{\substack{\nu:\Phi\to\lP^n\\ |\nu|=\al}}\prod_{f\in\Phi}
\psi_{\deg f}(r_{\nu(f)})(q).$$
We refer to Appendix for the basic definitions concerning
\la-rings. In the notation from Introduction, the formula of Kac
and Stanley  has the form

\begin{lem}\label{lmm Kac-St reformulation}
We have
$$m=\prod_{d\ge 1}\psi_d(r)^{\Phi_d}.$$
\end{lem}
\begin{proof}
Using the Kac-Stanley formula we get
\begin{align*}
\sum_{\al\in\cN^n}m_\al(q) x^\al
 =&\sum_{\nu:\Phi\to\lP^n}\prod_{f\in\Phi} \psi_{\deg
 f}(r_{\nu(f)})(q)x^{|\nu|}\\
 =&\sum_{\nu:\Phi\to\lP^n}\prod_{f\in\Phi} \psi_{\deg
 f}(r_{\nu(f)}x^{|\nu(f)|})(q)
 =\prod_{f\in\Phi}(\sum_{\la\in\lP^n}\psi_{\deg
 f}(r_{\la}x^{|\la|})(q))\\
=&\prod_{f\in\Phi}\psi_{\deg f}(r)(q)%
=\prod_{d\ge1}\psi_d(r)(q)^{\Phi_d(q)}.
\end{align*}
\end{proof}

\begin{lem}\label{lmm:rel a and m}
We have $m=\Exp(a)$.
%or, equivalently, $a(q)=\Log(m(q))$.
%=\sum_{n\ge1}\frac{\mu(n)}n\psi_n(\log(m(q))).$$
\end{lem}
\begin{proof}
Let $i_\al(q)$ denote the number of indecomposable representations
(do not confuse them with the absolutely indecomposable
representations) of dimension \al. Then it holds (see
\cite[p.91]{Kac1})
$$i_\al(q)=\sum_{d\ge 1}\sum_{k\mid d}\frac{\mu(k)}da_{\al/d}(q^{d/k}),$$
where the sum runs over those $d$ that divide the coefficients of
\al. On the other hand, it is clear that
$$\sum_\al m_\al(q)x^\al=\prod_\al(1+x^\al+x^{2\al}+\dots)^{i_\al(q)}%
=\prod_\al(1-x^\al)^{-i_\al(q)},$$
so we have to prove that
$$\sum_{\al}i_\al(q)\log\frac 1{1-x^\al}=\log(m(q))=\Psi(a(q)).$$
We have
\begin{align*}
&\sum_{\al}\sum_{r\ge1}i_\al(q)\frac{x^{r\al}}r
 =\sum_{\al}\sum_{r\ge1}\sum_{d\ge 1}\sum_{k\mid d}
 \frac{\mu(k)}da_{\al/d}(q^{d/k})\frac{x^{r\al}}r\\
 &\rEq^{{\be=\al/d}}_{{m=kr,n=d/k}}
 \sum_{\be}\sum_{m\ge1}\sum_{n\ge 1}
 a_{\be}(q^{n})\frac{x^{mn\be}}{mn}\sum_{k\mid m}
 {\mu(k)}%
 =\sum_{\be}\sum_{n\ge1}a_\be(q^n)\frac{x^{n\be}}n=\Psi(a(q)).
\end{align*}
\end{proof}

\begin{thm}[Hua's formula]\label{thr intro relations}
We have
$$a(q)=(q-1)\Log(r(q))$$
or, equivalently, $m(q)=\Pow(r(q),q-1)$.
\end{thm}
\begin{proof}
The minimal polynomial of a nonzero element of $\cF_{q^n}$ over $\cF_q$
is contained in $\Phi$ and has degree dividing $n$. Conversely, for any
$d\mid n$ and any $f\in\Phi$ of degree $d$ there are $d$ elements in
$\cF_{q^n}$ having $f$ as a minimal polynomial. This implies
$q^n-1=\sum_{d\mid n}d\Phi_d(q)$. Applying Lemma \ref{lmm:power formula}, we get
$$m(q)=\prod_{d\ge1}\psi_d(r(q))^{\Phi_d(q)}=\Pow(r(q),q-1).$$
\end{proof}

\begin{rem}
The proof of the above theorem can be generalized as follows 
(cf. \cite{Rod1}). Let $X$ be an algebraic variety over $\cF_q$ 
(i.e., a separated scheme of finite
type over $\cF_q$) such that there exists a polynomial $p_X$ satisfying
$\#X(\cF_{q^n})=p_X(q^n)$ for $n\ge1$. Let $|X|$ be the set of closed points
of $X$ and for any $x\in |X|$ let $\deg x=[k(x):\cF_q]$. Then
$\#\sets{x\in|X|}{\deg x=d}=\Phi_d(q)$, where the polynomials
$\Phi_d$, $d\ge1$ are defined by the formula
$\sum_{d\mid n}d\Phi_d=\psi_n(p_X)$. Given a function
$r\in \cQ[q]\pser{x_1,\dots,x_n}^+$, we have
$$\Pow(r,p_X)=\prod_{d\ge1}\psi_d(r)^{\Phi_d}.$$
Considering $q$ as a number of the elements of the base field, we get
$$\Pow(r,p_X)(q)=\prod_{x\in |X|}\psi_{\deg x}(r)(q).$$
\end{rem}

\section{On the verification of the Kac conjecture}\label{sec:criterium}
Let \g be a Kac-Moody algebra with a generalized Cartan matrix $C$
as in Introduction. Let \gh be a Cartan subalgebra of \g and let
$\set{\al_1,\dots,\al_n}\sb\gh^*$ be the simple roots of \g. Let
$(-,-)$ be the standard non-degenerate symmetric bilinear form on
$\gh^*$ (see \cite[Section 2.1]{Kac3}). Let $W$ be the Weyl group
of \g and $\rho\in \gh^*$ be any element with $(\rho,\al_i)=1$ for
all $i=1,\dots,n$. We always identify the root lattice (generated
by $\al_1,\dots,\al_n$) with $\cZ^n$. Note that the restriction of
the bilinear form $(-,-)$ to this lattice is given by the matrix
$C$. In particular, for any $\al\in\cN^n$ it holds $T(\al)=\frac
12(\al,\al)$.

The goal of this section is to prove Theorem \ref{thr intro
criterion}.
%\begin{theorem}\label{prp main in section2}
%The first Kac conjecture is true if and only if for any
%$\al\in\cN^n$ either $r_\al(0)=0$ or $T(\al)=\height\al$.
%\end{theorem}
One direction of the theorem (the ``only if'' condition) is rather
simple. Assume that $a_\al(0)$ coincides with the root
multiplicity $\mult\al$ for any $\al\in\cN^n$. Then we deduce from
Hua's formula that $r(q)=\Exp(\frac {a(q)}{q-1})$ and therefore
\begin{gather*}
r(0)=\Exp(-a(0))
=\prod_\al\Exp(x^\al)^{-a_\al(0)}\\
=\prod_\al(1+x^\al+x^{2\al}+\dots)^{-a_\al(0)}
=\prod_\al(1-x^\al)^{a_\al(0)}
=\sum_{w\in W}(-1)^{l(w)}x^{\rho-w\rho},
\end{gather*}
where the last equation is the Kac-Weyl denominator formula (see,
e.g., \cite[10.4.4]{Kac3}). We note that for any $\al\in\cN^n$ it
holds $T(\al)-\height\al=\frac 12(\al,\al)-(\rho,\al)$. In
particular,
\begin{gather}\notag
T(\rho-w\rho)-\height(\rho-w\rho)
=\frac 12 (\rho-w\rho,\rho-w\rho)-(\rho,\rho-w\rho)\\
\label{eq1}
%=\frac{-1}2(\rho+w\rho,\rho-w\rho) %
=\frac{1}2((w\rho,w\rho)-(\rho,\rho))=0,
\end{gather}
where the last equality follows from the $W$-invariance of the
bilinear form $(-,-)$. To prove the other direction we will apply
the approach of Peterson \cite{Pet1} for a recursive calculation
of the root multiplicities. Roughly speaking, one shows in this
approach that the generating function $\ub a$ of root
multiplicities (or rather the function $\Exp(-\ub a)$) satisfies
certain second order differential equation and can be determined
by its ``boundary values''. We will show that, under the
conditions of the proposition, the function $r(0)=\Exp(-a(0))$
also satisfies this differential equation. As it has the same
``boundary values'' as for the root multiplicities, we obtain
$a(0)=\ub a$.

The original paper by Peterson \cite{Pet1} is unpublished. Our
references for his approach were \cite[Exercise 11.11]{Kac3} and
\cite{Kang1}. First, we define some formal differential operators.
Define the derivation
$\Na:\cC\pser{x_1,\dots,x_n}\ar\gh^*\pser{x_1,\dots,x_n}$ by the
formula
$$\Na f:=\sum_\al f_\al\al x^\al,
\qquad \text{ for any }f=\sum_\al f_\al x^\al.$$ %
It is easy to see that it satisfies $\Na(fg)=f\Na(g)+\Na(g)f$.
Define the symmetric bilinear form
$$(-,-):\gh^*\pser{x_1,\dots,x_n}\xx \gh^*\pser{x_1,\dots,x_n}\ar
\cC^*\pser{x_1,\dots,x_n}$$ %
by the formula
$$(F,G):=\sum_{\al,\be}(F_\al,G_\be)x^{\al+\be},$$
where $F=\sum_\al F_\al x^\al$ and $G=\sum_\al G_\al x^\al$ are
elements in $\gh^*\pser{x_1,\dots,x_n}$. For any $\la\in\gh^*$,
define the operator $\dd_\la$ on $\cC\pser{x_1,\dots,x_n}$  by
$$\dd_\la f:=(\la,\Na f)=\sum_\al (\la,\al)f_\al x^\al.$$ %
\begin{sloppypar}
Let $\la_1,\dots,\la_r$ be an orthonormal basis of $\gh^*$. Define
the operator $\De$ on $\cC\pser{x_1,\dots,x_n}$ by the formula
\end{sloppypar}
$$\De f:=\sum_{i=1}^r \dd_{\la_i}^2f,\qquad f\in\cC\pser{x_1,\dots,x_n}.$$ %

%Given $f=\sum_\al f_\al x^\al$ and $g=\sum_\al g_\al x^\al$ in
%$\cC\pser{x_1,\dots,x_n}$, define
%$(f,g)\in\cC\pser{x_1,\dots,x_n}$ by the formula
%$$(f,g)=\sum_{\al,\be}(\al,\be)f_\al g_\be x^{\al+\be}.$$

\begin{lem}
Let $f=\sum_\al f_\al x^\al\in\cC\pser{x_1,\dots,x_n}$. Then
\begin{enumerate}
    \item $\De f=\sum_{i=1}^r \dd_{\la_i}^2f=\sum_\al (\al,\al)f_\al
    x^\al$.
    \item $\sum_{i=1}^r (\dd_{\la_i}f)^2=(\Na f,\Na f)$.
\end{enumerate}
\end{lem}
\begin{proof}
It is enough to prove the formula for $f=x^\al$. But
$$\De(x^\al)=\sum_{i=1}^r (\la_i,\al)^2x^\al=(\al,\al)x^\al$$
by the Pythagorean theorem. The second statement is analogous.
\end{proof}

\begin{lem}\label{lmm formulas for De and Na}
Let $f\in\cC\pser{x_1,\dots,x_n}^+$ and $g:=\exp(f)$. Then
\begin{enumerate}
    \item $\frac{\Na g}g=\Na f$,
    \item $\frac{\De g}g=\De f+(\Na f,\Na f)$.
\end{enumerate}
\end{lem}
\begin{proof}
From $f=\log(g)$ we obtain $\Na f=\Na\log(g)=\frac{\Na g}g$. Let
$\dd_i=\dd_{\la_i}$. Then $\dd_i\exp(f)=\exp(f)\dd_if$ and
therefore $\dd_i^2\exp(f)=\exp(f)\dd_i^2f+\exp(f)(\dd_i f)^2$. It
follows
$$\De\exp(f)=\sum_{i=1}^r\dd_i^2f=\exp (f)\sum_{i=1}^r(\dd_i^2 f+(\dd_i f)^2)
=\exp(f)(\De f+(\Na f,\Na f)).$$ 
\end{proof}

Let $\ub a_\al:=\mult \al$ be the multiplicities of the Kac-Moody
algebra \g. Define $\ub a:=\sum \ub a_\al x^\al$, $\ub c:=\Psi\ub
a$ and $\ub r:=\Exp(-\ub a)=\exp(-\ub c)$.

\begin{lem}
It holds
\begin{enumerate}
    \item $(\De-2\dd_\rho)\ub r=0$,
    \item $(\De-2\dd_\rho)\ub c=(\Na \ub c,\Na\ub c)$.
\end{enumerate}
\end{lem}
\begin{proof}
We note that
$$\ub r=\Exp(-\ub a)=\prod_{\al\in\cN^n}(1-x^\al)^{\ub a_\al}
=\sum_{w\in W}(-1)^{l(w)}x^{\rho-w\rho}.$$ %
To prove the first formula, we have to show that
$(\De-2\dd_\rho)(x^{\rho-w\rho})=0$, i.e.,
$$(\rho-w\rho,\rho-w\rho)-2(\rho,\rho-w\rho)=0.$$
But this was already shown in \eqref{eq1}. The second formula
follows from the first one if we recall that $\ub r=\exp(-\ub c)$
and apply Lemma \ref{lmm formulas for De and Na}. 
\end{proof}

This Lemma implies that for any $\al\in\cN^n$
$$(\al,\al-2\rho)\ub c_\al=\sum_{\be+\ga=\al}(\be,\ga)\ub c_\be \ub c_\ga.$$
We recall that $\ub c=\Psi\ub a$ and therefore $\ub
c_\al=\sum_{n\ge1} \frac 1n \ub a_{\al/n}$, where the sum is over
all $n$ dividing the coefficients of \al. We deduce that
$$(\al,\al-2\rho)\ub a_\al=\sum_{\be+\ga=\al}(\be,\ga)\ub c_\be \ub c_\ga
-(\al,\al-2\rho)\sum_{n>1} \frac 1n \ub a_{\al/n}$$ %
and this allows us to calculate the numbers $\ub a_\al$
inductively by height (note that $(\al,\al)<2(\rho,\al)$ for any
non-simple root by \cite[11.6.1]{Kac3}). The initial values (i.e.,
values at simple roots) are equal $1$.

Assuming that the condition
$$r_\al(0)=0\text{ or }T(\al)=\height\al,\qquad \forall\al\in\cN^n$$
of Theorem \ref{thr intro criterion} is satisfied, we will show
that the coefficients of $a(0)$ satisfy the same recursive formula
as the coefficients of $\ub a$ above. Let $c=\Psi(a)$.

\begin{lem}\label{lmm formulas for r and c}
Assume that for any $\al\in\cN^n$ either $r_\al(0)=0$ or
$T(\al)=\height\al$. Then
\begin{enumerate}
    \item $(\De-2\dd_\rho)r(0)=0$.
    \item $(\De-2\dd_\rho)c(0)=(\Na c(0),\Na c(0)).$
\end{enumerate}
\end{lem}
\begin{proof}
To prove the first formula, we have to show that for any
$\al\in\cN^n$ it holds
$$(\De-2\dd_\rho)(r_\al(0)x^\al)=((\al,\al)-2(\rho,\al)))r_\al(0)x^\al=0,$$
but this is precisely the condition of the lemma. The second
formula follows from the first one by applying Lemma \ref{lmm
formulas for De and Na}. 
\end{proof}

\begin{proof}{Proof of Theorem \ref{thr intro criterion}.}
The ``only if'' part has already been shown in the beginning of
this section. Assume that for any $\al\in\cN^n$ either
$r_\al(0)=0$ or $T(\al)=\height\al$. Then it follows from Lemma
\ref{lmm formulas for r and c} that
$$(\De-2\dd_\rho)c(0)=(\Na c(0),\Na c(0))$$
and we obtain, in the same way as above, that the coefficients of
$a(0)$ can be inductively determined by the formula
$$(\al,\al-2\rho)a_\al(0)=\sum_{\be+\ga=\al}(\be,\ga)c_\be(0)
c_\ga(0)-(\al,\al-2\rho)\sum_{n>1} \frac 1n a_{\al/n}(0).$$ %
Here we note that if $\al\in\cN^n$ is not a root then $a_\al=0$ by
the results of Kac \cite{Kac1} and if \al is a non-simple root
then $(\al,\al)<2(\rho,\al)$ as it has already been mentioned
above. The initial values (i.e., values at simple roots) equal $1$
as before. All this implies that $a(0)=\ub a$ and the theorem is
proved. 
\end{proof}

\begin{exmp}\label{exm test}
Consider a quiver with the underlying graph
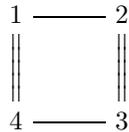
\begin{diagram}
1&\rLine&2\\
\dEq&&\dEq\\
4&\rLine&3
\end{diagram}
Using the formula for $r_\al(q)$ given in Introduction, we can write the first
terms of $r(0)$ (coefficients by $x^\al$ with $\al$ smaller than $(3,3,3,3)$) 
\begin{align*}
r(0)
&=x^{(0,0,0,0)}-x^{(0,0,0,1)}-x^{(0,0,1,0)}-x^{(0,1,0,0)}
-x^{(1,0,0,0)}+x^{(0,1,0,1)}+x^{(1,0,1,0)}\\
&+x^{(0,0,1,2)}+x^{(0,0,2,1)}+x^{(1,2,0,0)}+x^{(2,1,0,0)}
-x^{(0,0,2,2)}-x^{(2,2,0,0)}+x^{(0,1,3,0)}\\
&+x^{(0,3,1,0)}+x^{(1,0,0,3)}+x^{(3,0,0,1)}
-x^{(0,3,1,2)}-x^{(1,2,0,3)}-x^{(2,1,3,0)}-x^{(3,0,2,1)}.
\end{align*}
It is easy to see that all $\al$ with $r_\al(0)\ne0$ satisfy $T(\al)=\height(\al)$.
\end{exmp}

\appendix
\section{\la-Rings}
We follow Getzler \cite{Getz1}. Let \La be the ring of
symmetric functions, $e_n$ be the elementary symmetric functions,
$h_n$ be the complete symmetric functions, and $p_n$ be the power
sums (see, e.g., \cite{Mac1}).

Define the operation $\circ:\La\xx\La\ar\La$, called plethysm, by
the properties
$$p_n\circ f(x_1,\dots x_k)=f(x_1^n,\dots,x_k^n), \quad n\ge1,\ f\in\La,$$
$$(-\circ f):\La\ar\La\text{ is a ring homomorphism for any }f\in\La.$$

\begin{rem}
Note that this operation is associative. It is, however, not
commutative and not additive in the second argument. For example,
$0\circ 1=0$ and $1\circ 0=1$.
\end{rem}

\begin{defn}
A pre-$\la$-ring is a commutative ring $R$ together with a map
${\circ:\La\xx R\ar R}$ such that $(-\circ a):\La\ar R$ is a ring
homomorphism for any $a\in R$ and, for $\la_n:=(e_n\circ -):R\ar
R$, it holds $\la_1=\Id_R$ and
$\la_n(a+b)=\sum_{i=0}^n\la_i(a)\la_{n-i}(b)$.
\end{defn}

\begin{rem}
The operations $\la_n:R\ar R$, $\si_n:R\ar R$ and $\psi_n:R\ar R$
induced, respectively, by $e_n$, $h_n$ and $p_n$ are called
\la-operations, \si-operations and Adams operations, respectively.
To define the pre-\la-ring structure on $R$ it is enough just to
define the \la-operations or \si-operations satisfying the
corresponding conditions. If $R$ is an algebra over \cQ then it
suffices to define the Adams operations satisfying $\psi_1=\Id_R$
and $\psi_n(a+b)=\psi_n(a)+\psi_n(b)$. For simplicity, we will
always assume that $R$ is an algebra over \cQ.
\end{rem}

\begin{defn}
A pre-\la-ring $R$ is called a \la-ring if it holds
$$f\circ(g\circ a)=(f\circ g)\circ a,\quad f,\ g\in\La,\ a\in R$$
and
%the multiplication map $R\ts R\ar R$ is a map of pre-\la-rings.
$\la_n(1)=0$ for any $n\ge2$.
\end{defn}

\begin{rem}
A pre-\la-ring $R$ is a \la-ring if and only if $\psi_n(1)=1$ for
every $n\ge1$ and $\psi_m(\psi_n(a))=\psi_{mn}(a)$ for any $m,n\ge
1$ and $a\in R$. It can be shown that $\psi_n$ are actually 
ring homomorphisms.
\end{rem}

\begin{exmp}
The basic example is a \la-ring with all Adams operations being
identities. If $R$ is a \la-ring, we endow the ring
$R[x_1,\dots,x_r]$ with a \la-ring structure by defining the Adams
operations as
$$\psi_n(a\cdot x^\al)=\psi_n(a)\cdot x^{n\al},\quad n\ge1,\ a\in R,\ \al\in\cN^r.$$
In the same way, we endow the ring of
formal power series over $R$ with a \la-ring structure.
\end{exmp}

\begin{rem}
In order to work with infinite sums in the \la-rings, we will
assume that they are complete graded \la-rings. A \la-ring $R$ is
called a complete graded \la-ring if it is a complete graded ring
$R=\hat\oplus_{n\ge0}R_n$ and $\psi_m(R_n)\sb R_{mn}$ for any
$m\ge 1$, $n\ge 0$. We define $R^+=\hat\oplus_{n\ge1}R_n$. Our
main example of a complete graded \la-ring is the ring
$R\pser{x_1,\dots,x_r}$, where $R$ is a usual \la-ring and the
grading is given by
$$\deg(ax^\al)=\height\al=\sum \al_i,\quad a\in R,\ \al\in\cN^r.$$
\end{rem}

\begin{lem}
Let $R$ be a complete graded \la-ring. The map $\Psi:R^+\ar R^+$,
$\Psi(f)=\sum_{n\ge1}\frac 1n\psi_n(f)$ has an inverse
$\Psi\inv:R^+\ar R^+$ given by $\Psi\inv(f)=\sum_{n\ge1}\frac
{\mu(n)}n\psi_n(f),$ where $\mu(n)$ is a M{\"o}bius function.
\end{lem}
\begin{proof}
We will just show that $\Psi\inv\Psi=\Id$, as the equality
$\Psi\Psi\inv=\Id$ is analogous. The basic property of the
M{\"obius} function is that $\sum_{k\mid n}\mu(k)=0$ for $n\ne 1$.
We deduce
$$\Psi\inv\Psi(f)=\sum_{k,m\ge1}\frac{\mu(k)}k\frac 1m\psi_{km}(f)=
\sum_{n\ge1}\frac{\psi_n(f)}n\sum_{k\mid n}\mu(k)=\psi_1(f)=f.$$
\end{proof}

In the next corollary, we use the maps $\exp:R^+\ar 1+R^+$ and
$\log:1+R^+\ar R^+$. Their definition can be found, e.g., in
\cite[Ch.II \S 6]{BourLie13}. We define the map ${\Exp:R^+\ar
1+R^+}$ by
$$\Exp(f):=\sum_{n\ge0}\si_n(f),\qquad f\in R^+.$$
It is easy to see that $\Exp(f+g)=\Exp(f)\Exp(g)$. 
One knows that (see, e.g., \cite[2.10]{Mac1})
$$\sum_{k\ge 1}p_kt^{k-1}=\frac d{dt}\log \sum_{k\ge0}h_kt^k.$$ 
This implies
$\sum_{k\ge0}h_kt^k=\exp(\sum_{k\ge 1}\frac{p_kt^k}k)$ and
therefore
$$\Exp(f)=\exp\Big(\sum_{k\ge 1}\frac{\psi_k(f)}k\Big)=\exp(\Psi(f)).$$

\begin{cor}[Cadogan formula, see \cite{Cadogan1,Getz1}]
Let $R$ be a complete graded \la-ring. Then the map $\Exp:R^+\ar
1+R^+$ has an inverse $\Log:1+R^+\ar R^+$,
$\Log(f)=\Psi\inv(\log(f))$.
\end{cor}

Define the map $\Pow:(1+R^+)\xx R\ar 1+R^+$ (called a power structure in \cite{GLM2})
by the formula 
$$\Pow(f,g):=\Exp(g\Log(f)).$$ 
Analogously, define $f^g:=\exp(g\log(f))$.

\begin{lem}\label{lmm:power formula}
Let $f\in 1+R^+$, $g\in R$. Define the elements $g_d\in R$, $d\ge1$ by the formula 
$\sum_{d\mid n}d\cdot g_d=\psi_n(g)$. Then we have
$$\Pow(f,g)=\prod_{d\ge1}\psi_d(f)^{g_d}.$$
\end{lem}
\begin{proof}
After taking logarithms we have to prove
$$\Psi(g\Log(f))=\sum_{d\ge1}g_d\psi_d(\log(f)).$$
Let $h=\Log(f)$. Then $\Psi(h)=\log(f)$ and we have to show
$$\Psi(gh)=\sum_{d\ge1}g_d\psi_d\Psi(h).$$
We have $\psi_n(gh)=\psi_n(g)\psi_n(h)=\sum_{d\mid n}d g_d \psi_n(h)$.
This implies 
$$\Psi(gh)=\sum_{d\mid n}d g_d \frac{\psi_n(h)}n
=\sum_{d,k\ge1}g_d\frac{\psi_d\psi_k(h)}k=\sum_{d\ge1}g_d\psi_d\Psi(h).$$
\end{proof}

\bibliography{fullbib}
%\bibliographystyle{../tex/hamsplain}
%\bibliography{fullbib}
\bibliographystyle{hamsplain}
%GATHER{../tex/fullbib.bib}

\end{document}